\documentclass[11pt]{article}
\input amssym.def
\usepackage{graphicx}
\usepackage[latin1]{inputenc}

\def\sqw{\hbox{\rlap{\leavevmode\raise.3ex\hbox{$\sqcap$}}$%
\sqcup$}}

\def\codim{{\rm codim}\,}

  \def\bn{\hbox{\it 
I\hskip -2pt N}}      
\def\demo{\noindent{\bf Proof \ }} \newtheorem{theorem}{Theorem} 
\newtheorem{lemma}{Lemma} \newtheorem{proposition}{Proposition} 
\newtheorem{definition}{Definition} \newtheorem{example}{Example} 
  
\newtheorem{corollary}{Corollary}

\def\codim{{\rm \ codim}\,}
\def\dim{{\rm \ dim}\,}
\def\projdim{{\rm \ projdim}\,}

\def\cd{{\rm \ cd} }

\def\hut{{\rm \ ht} }
\def\pd{{\rm \ projdim}\  }
\def\depth{{\rm \ depth}\, }

\def\ara{{\rm \ ara\,} }
\def\mod{{\rm \ mod} }

\def\reg{{\rm \ reg}\, }

\def\bn{\hbox{\it I\hskip -2pt N}}

\def\ux{\underline x}

\newcommand{\Ical}{{\mathcal I}}
\newcommand{\Jcal}{{\mathcal J}}
\newcommand{\Bcal}{{\mathcal B}}
\newcommand{\Acal}{{\mathcal A}}
\newcommand{\Kcal}{{\mathcal K}}

\newcommand{\Pcal}{{\mathcal P}}

\newcommand{\Qcal}{{\mathcal Q}}

\newcommand{\Dcal}{{\mathcal D}}

\def\demo{\noindent{\bf Proof \ }} 

\begin{document}

\begin{center}
\uppercase{{\bf  $p-$Ferrer diagram, $p-$linear ideals and arithmetical rank }}
\end{center}
\advance\baselineskip-3pt
\vspace{2\baselineskip}
\begin{center}
{{\sc Marcel
Morales}\\
{\small Universit\'e de Grenoble I, Institut Fourier, 
UMR 5582, B.P.74,\\
38402 Saint-Martin D'H\`eres Cedex,\\
and IUFM de Lyon, 5 rue Anselme,\\ 69317 Lyon Cedex (FRANCE)}\\
 }
\vspace{\baselineskip}
\end{center}
{\bf Abstract} 

In this paper we introduce $p-$Ferrer diagram, note that $1-$ Ferrer diagram are the usual Ferrer diagrams or Ferrer board, and corresponds to planar partitions. To any $p-$Ferrer diagram we associate a $p-$Ferrer ideal. We prove that $p-$Ferrer ideal have Castelnuovo mumford regularity $p+1$. We also study Betti numbers , minimal resolutions of $p-$Ferrer ideals. Every $p-$Ferrer ideal is $p-$joined ideals  in a sense defined in a fortcoming paper \cite{m2}, which extends the notion of linearly joined ideals introduced and developped in the papers \cite{bm2}, \cite{bm4},\cite{eghp} and \cite{m1}.
We can observe the connection between the results on this paper about the Poincaré series of a $p-$Ferrer diagram $\Phi $and the rook problem, which consist to put $k$ rooks in a non attacking position on the $p-$Ferrer diagram $\Phi $. 

\section{Introduction}
We recall that any non trivial ideal $\Ical\subset S$ has a finite free resolution :
$$0\rightarrow F_s\buildrel{M_s}\over\rightarrow F_{s-1}\rightarrow ....\rightarrow F_1\buildrel{M_1}\over\rightarrow  \Ical\rightarrow 0$$ the number $s$ is called the projective dimension of $S/\Ical$ and the Betti numbers are defined by
$\beta_i(\Ical)=\beta_{i+1}(S/\Ical)= rank F_{i+1}$. By the theorem of Auslander and Buchsbaum we know that $s=\dim S-\depth(S/\Ical)$.
We will say that the ideal $\Ical$  has a pure resolution if  $F_{i}=S^{\beta _i}(-a_i)$ for all $i=1,...,s$.  This means that $\Ical$ is generated by elements in degree $a_1$, and for $i\geq 2$  the matrices $M_i$ in the minimal free resolution of $\Ical$ have homogeneous entries of degree $a_{i}-a_{i-1}$.

 We will say that the ideal $\Ical$  has a $p-$linear resolution if its minimal free resolution is linear,  i.e.  $\Ical$  has a pure resolution and  for $i\geq 2$  the matrices $M_i$ have linear entries.
 
If $\Ical$  has a pure resolution, then the  Hilbert series of $S/\Ical$ is given by:
$$H_{S/\Ical}(t)=\displaystyle\frac{1-\beta_1 t^{a_1}+...+ (-1)^s \beta_s t^{a_s}}{(1-t)^n}$$
where $n=\dim S$. Since $a_1<...<a_s$ it the follows that if $\Ical$  has a pure resolution then the Betti numbers are determined by the Hilbert series.

{\bf $p-$Ferrer partitions and  diagrams.}
The $1-$Ferrer partition is a nonzero natural integer $\lambda $, a $2-$Ferrer partition is called a partition and is given by a sequence $\lambda _1\geq ...\geq \lambda _m>0$ of natural integers, a $3-$Ferrer partition is called  planar partition. $p-$Ferrer partitions are defined inductively  $\Phi: \lambda_1 \geq  \lambda_2 \geq   ...\geq  \lambda_m$, where $\lambda_j$ is a 
$p-1$Ferrer partition for $j=1,...,m$, and the relation  $\leq $ is also defined recursively: if  $\lambda_i:\lambda_{i,1}\geq ...\geq \lambda_{i,s},
\lambda_{i+1}:\lambda_{i+1,1}\geq ...\geq \lambda_{i+1,s'}$  we will say that $\lambda_i\geq \lambda_{i+1}$
 if and only if $s\geq s'$ and $\lambda_{i,j}\geq \lambda_{i+1,j}$ for any $j=1,...,s'$.
Up to my knowledge there are very few results for $p-$Ferrer partitions in bigger dimensions.

To any $p-$Ferrer partition we associate a $p-$Ferrer diagram which are subsets of $\bn^p$. The 1-Ferrer diagram associated to $\lambda \in \bn$ is the subset $\{1,...,\lambda \}$. 
Inductively  if  $\Phi: \lambda_1 \geq  \lambda_2 \geq   ...\geq  \lambda_m$, is a $p-$Ferrer partition, where $\lambda_j$ is a 
$p-1$-Ferrer partition  for $j=1,...,m$, we associate to $\Phi $ the  $p-$Ferrer diagram
$\Phi =\{(\eta,1), \eta \in \lambda_1\}\cup ...\cup \{(\eta,m ), \eta \in \lambda_m\}$. Ferrer $p-$diagrams can also be represented by a set of boxes  labelled by a $p-$uple  $(i_1,...,i_p)$ of non zero natural  numbers, they have the property that if $1\leq i'_1\leq i_1,..., 1\leq i'_p\leq i_p$, then the 
box labelled $(i'_1,...,i'_p)$ is also in the $p-$Ferrer diagram. We can see that for two Ferrer diagrams: 
$\Phi _1 \geq \Phi _2$ if and only if the set of boxes of $\Phi _1$ contains the set of boxes of $\Phi _2$.

\begin{example}\label{exo4322}
The following picture corresponds to  the $3-$Ferrer diagram given by: 
 $$\begin{array}{llll}
4 & 3& 2& 2  \\
3 & 2& 1& 0  \\
2 & 0& 0& 0  \\
2 & 0& 0& 0  \\
\end{array}$$
 \begin{center}
\includegraphics[height=2 in,width=2in]{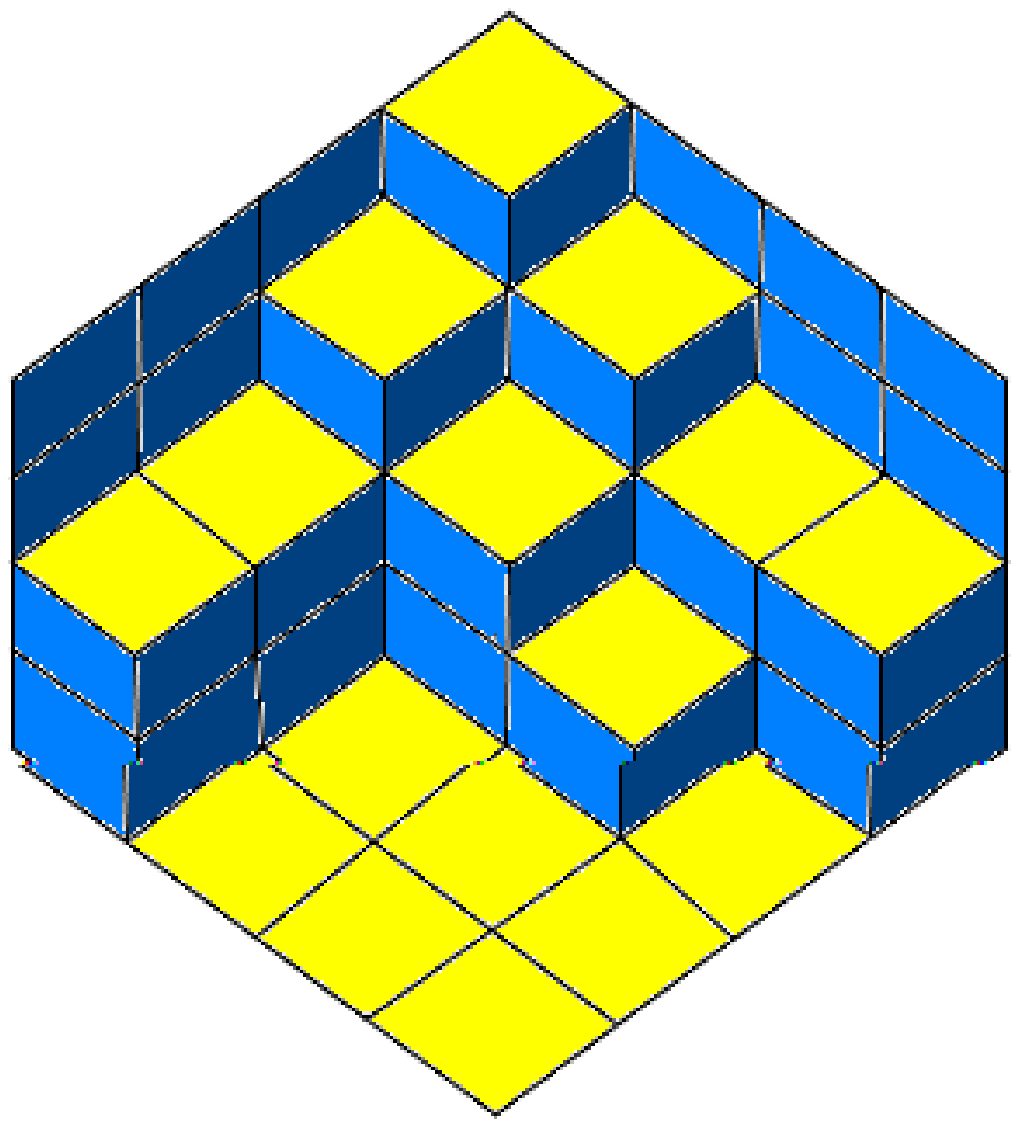}
\end{center}
\end{example}

 \begin{example}\label{exo54432}The following picture corresponds to  the $3-$Ferrer diagram given by: 
 $$\begin{array}{lllll}
5 & 4& 4& 3 & 2  \\
4 & 4& 3& 3 & 1  \\
4 & 4& 3& 1 & 0  \\
2 & 1& 1& 0 & 0  \\
2 & 1& 0& 0 & 0  \\
\end{array}$$
 \begin{center}
\includegraphics[height=2 in,width=2in]{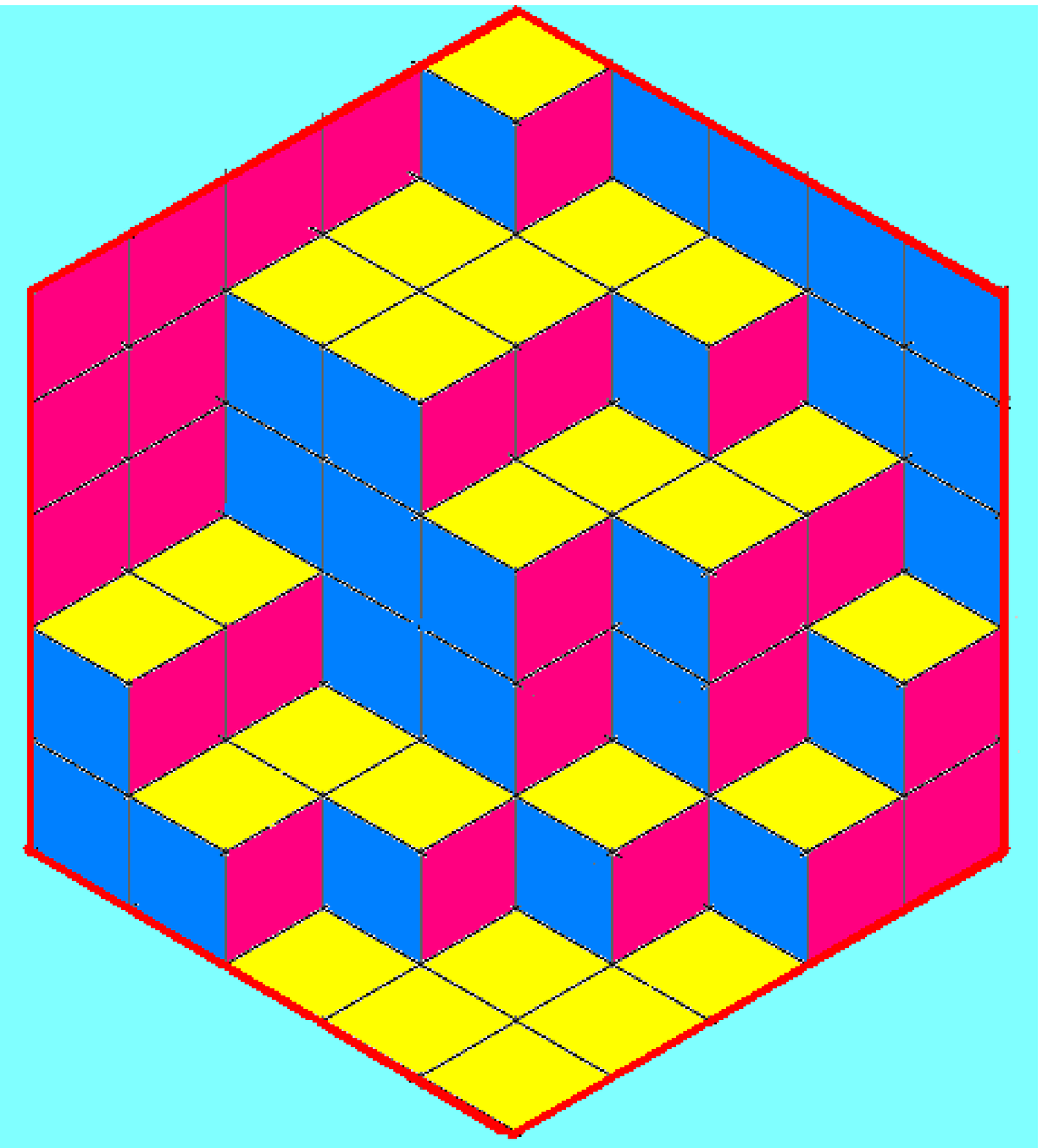}
\end{center}
\end{example}
\begin{definition} Given a $p-$Ferrer diagram (or partition) $\Phi$   we can associated a monomial ideal $\Ical_\Phi $ in the following way.
 Let  consider the polynomial ring $K[\ux^{(1)},\ux^{(2)},...,\ux^{(p)}]$ where $\ux^{(i)}$ stands for 
the infinitely set of variables :$\ux^{(i)}=  \{x_1^{(i)},x_2^{(i)},...   \},$ we define inductively the ideal $\Ical_\Phi $
\begin{enumerate}
\item For $q=2$ let $\Phi :\lambda\in \bn^*$, then $\Ical_\Phi $  is the ideal generated by the variables 
$ x_1^{(1)},x_2^{(1)},...,x_\lambda^{(1)}.$
\item For $q=2$ let $\lambda_1 \geq  \lambda_2 \geq   ...\geq  \lambda_m$ be a $2-$Ferrer diagram, then $\Ical_\Phi $ is an ideal
in the ring of polynomials $K[x_1,...,x_m, y_1,...,y_{\lambda_1 }]$
 generated by the monomials $x_i y_j$ such that $i=1,...,m$ and $j=1,...,\lambda_i$. In this case $x_j^{(1)}=y_j, x_j^{(2)}=x_i. $

\item For $p>2$ let $\lambda_1 \geq  \lambda_2 \geq   ...\geq  \lambda_m$ be a Ferrer diagram, where $\lambda_j$ is a 
$p-1$Ferrer diagram. Let $\Ical_{\lambda_j} \subset K[\Lambda ]$ be the ideal associated to $\lambda_j$, where $K[\Lambda ]$ 
s a polynomial ring in a finite set of variables  then $\Ical_\Phi $ is an ideal
in the ring of polynomials $K[x_1^{(p)},...,x_m^{(p)}, \Lambda ]$
 generated by the monomials $x_i^{(p)} y_j$ such that $i=1,...,m$ and  $y_j\in \Ical_{\lambda_i} $. That is
$$\Ical_{\Phi }= (\bigcup_{i=1}^m \{x_{i}^{(p)}\}\times \Ical_{\lambda_i }).$$
\end{enumerate}

\end{definition}
We can observe the connection between the results on this paper about the Poincaré series of a $p-$Ferrer diagram $\Phi $and the rook problem, which consist to put $k$ rooks in a non attacking position on the $p-$Ferrer diagram $\Phi $. This will be developped in a forthcoming paper.
 \section{$p-$Ferrer' ideals}
\begin{lemma} Let $S$ be a polynomial ring, $\Gamma_2...,\Gamma_r$be  non empty disjoint sets of variables, set $\Acal_i$ the ideal generated by $\Gamma_{i+1},..., \Gamma_r$. Let $\Bcal_2\subset ...\subset \Bcal_r$ be a sequence of ideals (not necessarily distinct), generated by the  sets $B_2\subset ...\subset B_r$. We assume that  no variable of  $\Gamma_{2}\cup ...\cup  \Gamma_r$  appears in $B_2, ..., B_r$, then 
$$ \Acal_1\cap (\Acal_2,\Bcal_2)\cap ... \cap (\Bcal_r)= (\bigcup_{2}^r \Gamma_i\times B_i )$$
where for two subsets $A,B\subset S$, we have set $A\times B=\{a\ b\ \mid \ a\in A,\ b\in B\}.$
\end{lemma}
\demo Let remark that if $\Gamma $ is a set of variables and $P\subset S$ is a set of polynomials such that no variable of $\Gamma $  appears in the elements of $P$ then $(\Gamma )\cap (P)=(\Gamma \times P)$. Moreover  if  $\Gamma_1, \Gamma _2 $ are disjoint sets of variables  and $P\subset S$ is a set of polynomials such that no variable of $\Gamma_1, \Gamma _2 $   appears in the elements of $P$ then $(\Gamma_1, \Gamma _2)\cap (\Gamma_1, ,P)=(\Gamma_1, \Gamma_2 \times P).$

We prove  by induction on the number $k$ the following statement:
$$ \Acal_1\cap (\Acal_2,\Bcal_2)\cap ... \cap (\Acal_k,\Bcal_k)= (\Acal_k,\bigcup_{2}^k \Gamma_i\times B_i ).$$ 
If $k=2$, it is clear that $\Gamma _2\times B_2 \subset \Acal_1\cap (\Acal_2,\Bcal_2)$, now let $f\in \Acal_1\cap (\Acal_2,\Bcal_2)$, we can write $f=f_1+f_2$, where $f_1\in (A_2)$,  $f_2\in (\Gamma _2)$ and no variable of $A_2$ appears in $f_2$, it follows that $f_2\in (\Gamma _2)\cap (B_2)=(\Gamma _2\times B_2 )$. 

Suppose that $$ \Acal_1\cap (\Acal_2,\Bcal_2)\cap ... \cap (\Acal_k,\Bcal_k)= (\Acal_k,\bigcup_{2}^k \Gamma_i\times B_i ),$$ we will prove that 
$$ \Acal_1\cap (\Acal_2,\Bcal_2)\cap ... \cap (\Acal_{k+1},\Bcal_{k+1})= (\Acal_{k+1},\bigcup_{2}^{k+1} \Gamma_i\times B_i ).$$ 
 Since $\Gamma_i\subset \Acal_j,$ for $j<i$, and $B_i\subset B_j$ for $i\leq j$,
   we have $\bigcup_{2}^{k+1} \Gamma_i\times B_i\subset (\Acal_j,\Bcal_j)$ for $1\leq j\leq k$, so we have the inclusion $"\supset "$.
   
By induction hypothesis we have that 
$$ \Acal_1\cap (\Acal_2,\Bcal_2)\cap ... \cap (\Acal_{k+1},\Bcal_{k+1})=(\Acal_k,\displaystyle\bigcup_{2}^k \Gamma_i\times B_i ) \cap (\Acal_{k+1},\Bcal_{k+1}).$$

Now let $f\in (\Acal_k,\displaystyle\bigcup_{2}^k \Gamma_i\times B_i ) \cap (\Acal_{k+1},\Bcal_{k+1}).$ we can write $f=f_1+f_2+f_3$, where $f_3\in (\bigcup_{2}^k \Gamma_i\times B_i) \subset \Bcal_{k+1}$, $f_1\in \Acal_{k+1}$, and $f_2\in (\Gamma_{k+1})$, and no variable of  $\Gamma_{k+1}\cup ...\cup \Gamma _r$ appears in $f_2$, this would imply that $f_2\in (\Gamma_{k+1})\cap \Bcal_{k+1}=(\Gamma_{k+1}\times B_{k+1}).$

\begin{definition}Let $\lambda_{m+1}=0,\delta_0=0,$ $\delta_1$ be the highest integer such that $\lambda_1=...=\lambda_{\delta_1}, $
 and by induction we define 
$\delta_{i+1}$ as the highest integer such that $\lambda_{\delta_i+1}=...=\lambda_{\delta_{i+1}},$ and set $l$ such that
$\delta_{l-1}=m.$
For $i=0,...,l-2$ let
$$\Delta_{l-i}=\{x_{\delta_i+1}^{(p)} ,... , x_{\delta_{i+1}}^{(p)}  \}, \Pcal_{l-i}=
\Ical_{\lambda _{\delta _{i+1}}}.$$
\end{definition}
So we have:
$\Phi =\{(\eta,1), \eta \in \lambda_1\}\cup ...\cup \{(\eta,m ), \eta \in \lambda_m\}$ and
$$\Ical_\Phi =(\bigcup_{i=2}^l \Delta_{i}\times P_{i} )=(\bigcup_{i=1}^m \{x_{i}^{(p)}\}\times \Ical_{\lambda_i }).$$ 
where for all $i$,  $P_i$ is a set of generators of $\Pcal_{i}$.

The  following Proposition is an immediate consequence of the above lemma :
\begin{proposition} 
\begin{enumerate}
\item We have the following decomposition (probably redundant):
$$  \Ical_{\Phi}=(x_{1}^{(p)},...,x_{m}^{(p)})\cap (x_{1}^{(p)},...,x_{m-1}^{(p)},\Ical_{\lambda_m})...
\cap (x_{1}^{(p)},...,x_{i-1}^{(p)},\Ical_{\lambda_i})\cap... (\Ical_{\lambda_m}),$$
\item Let $\Dcal_i=(\bigcup_{j=i+1}^{l}\Delta_{j})$, and $\Qcal_i= ( \Dcal_i,\Pcal_{i}).$ Then $$ \Ical_\Phi =\Qcal_1\cap \Qcal_2\cap ...\cap \Qcal_l.$$
\item  The minimal primary decomposition of $\Ical_\Phi $ is obtained inductively. Let $\Ical_{\lambda_{\delta_{i}}}=
\Qcal_1^{(i)}\cap ...\cap \Qcal_{r_i}^{(i)} $ be a minimal prime decomposition, where by induction hypothesis 
$\Qcal_{j}^{(i)}$ is a linear
ideal, then  the minimal  prime decomposition of $\Ical_\Phi$ is obtained from this decomposition by putting out unnecessary components.
\end{enumerate}
\end{proposition}
\begin{example} let $\Pcal_2=(c,d)\cap (e)$, $\Pcal_3=(c,d)\cap (c,e)\cap (e,f)$ and
$$\Ical_\Phi=(a,b)\cap (a,\Pcal_2)\cap \Pcal_3 $$ then 
$$\Ical_\Phi=(a,b)\cap (a,e)\cap (c,d)\cap (c,e)\cap (e,f) $$ is its minimal prime decomposition.
\end{example}

\begin{proposition} Let $\Ical\subset R$ be a $p-$Ferrer ideal then $\reg(\Ical)=p=\reg(R/\Ical)+1$.
\end{proposition}
\demo
 For any two ideals $\Jcal_1,\Jcal_2\subset S$ we have the following exact sequence:
$$0\rightarrow S/\Jcal_1 \cap \Jcal_2\rightarrow S/\Jcal_1\oplus  S/\Jcal_2\rightarrow S/(\Jcal_1 + \Jcal_2)\rightarrow 0$$
 From \cite[p. 289]{br-sh} 
$$\reg (S/\Jcal_1 \cap \Jcal_2)\leq \max \{ \reg (S/\Jcal_1\oplus  S/\Jcal_2), \reg(S/(\Jcal_1 + \Jcal_2))+1\}$$
 in our case we take $\Jcal_1=\bigcap_{i=1}^k  \Qcal_i, \Jcal_2=\Qcal_{k+1}$, so that
 $\reg(S/(\bigcap_{i=1}^k  \Qcal_i + \Qcal_{k+1}))=\reg(S/(\Dcal_{k}+\Pcal_{k+1}))=\reg(S'/(\Pcal_{k+1}))=p-1$, where $S=S'[\Dcal_{k}]$.
 It then follows that $\reg(S/(\bigcap_{i=1}^l  \Qcal_i)) \leq p$, on the other hand  $(\bigcap_{i=1}^l  \Qcal_i))$ is generated by elements
 of degree $p$, this implies $\reg(S/(\bigcap_{i=1}^l  \Qcal_i)) = p$.

We will show that in fact $\pd (S/\Ical_\lambda)$ is the number of diagonals in a $p-$Ferrer diagram.
\begin{definition}Let $\Phi  :\lambda_1 \geq  \lambda_2 \geq   ...\geq  \lambda_m$ be a $p-$Ferrer diagram.
We will say that the monomial in the $p-$Ferrer ideal (or diagram) $x_{\alpha _{p}}^{(p)}x_{\alpha _{p-1}}^{(p-1)}...x_{\alpha _{1}}^{(1)} $ is in the 
$\alpha _{p}+\alpha _{p-1}+...+\alpha _{1}-p+1$ diagonal.
Let $s_\Phi (k)$ be the number of elements in the $k-$diagonal of $\Phi $,
we will say that the $k-$diagonal of $\Phi $ is full if $s_\Phi (k)={{k-1+p-1}\choose{p-1}}$, which is the number of elements in the $k-$diagonal of $\bn^p$, let remark that by the definition of $p-$Ferrer diagram if the $k-$diagonal of $\Phi $ is full then the $j-$diagonal of $\Phi $ is full  for all $j=1,...,k$.
\end{definition}
\begin{lemma}
\begin{enumerate}
\item We have the formula 
$$ s_\Phi (k)=\sum_{i=1}^m s_{\lambda_i }(k-(i-1)), $$
\item Let $df(\Phi ) $ be the number of full diagonals of $\Phi $, then 
$$df(\Phi ) = \min \{ df(\lambda _i)+i-1 \ \mid \ i=1,...,m\}$$
\item Let $\delta (\Phi ) $ be the number of  diagonals of $\Phi $, then 
$$\delta (\Phi ) = \max \{ \delta (\lambda _i)+i-1 \ \mid \ i=1,...,m\},$$ and 
$\delta  (\Phi )=\max_{i=2}^l \{\delta  (\Pcal_i)+ \dim\Dcal_{i-1} -1  \}.$
\end{enumerate}

\end{lemma}
\demo The first item counts the number of elements in the $k-$diagonal of $\Phi$ by counting  all the $i-$ slice pieces. The second item means that the $k-$diagonal of $\Phi$ is full if and only if the $k-(i-1)-$diagonal of the $i-$ slice piece is full, and finally the third item means there is an element in the $k-$diagonal of $\Phi$ if and only if there is at least one element in the $k-(i-1)-$diagonal of the $i-$slide piece of $\Phi $, for some $i$. 

Remark that  
$\delta  (\Phi )=\max_{i=2}^l \{\delta  (\Pcal_i)+ \dim\Dcal_{i-1} -1  \},$
since  $\max_{i=1}^{\delta_1 } \{\delta  (\lambda_i)+i-1  \}=\delta  (\lambda_1)+\delta_1-1=
\delta(\Pcal_l)+\dim \Dcal_{l-1}-1$,
$\max_{i=\delta_1+1}^{\delta_2 } \{\delta  (\lambda_i)+i-1  \}=\delta  (\lambda_{\delta_1+1})+\delta_1+\delta_2-1=
\delta(\Pcal_{l-1})+\dim \Dcal_{l-2}-1$, and so on.
\begin{theorem}
Let consider a $p-$Ferrer diagram $\Phi$ and its associated ideal   $\Ical_\Phi $ in a polynomial ring $S$.
 Let $n=\dim S, c=\hut \Ical_\Phi, $ $d=n-c$. For $i=1,...,d-\depth S/\Ical$,  let $s_{d-i}$ be the numbers of elements in the
 $c+i$ diagonal of $\Phi $. 
Then :\begin{enumerate} 
\item $c$ the height of $\Ical_\Phi $ is equal to the number of full diagonals.
\item For $j \geq 1 $ we have $$\beta_j(S/\Ical_\Phi)= 
{{c+p-1}\choose{j+p-1}} {{j+p-2}\choose{p-1}}+\sum_{i=0}^{d-1}s_i {{n-i-1}\choose{j-1}}$$
\item $\pd (S/\Ical_\Phi )= \delta  (\Phi ) .$
 \end{enumerate}
 
\end{theorem}
\demo 
 \begin{enumerate}
\item We prove the statement by induction on $p$, if $p=1$ and $\Phi =\lambda \in \bn$, then  $\Ical_\phi=(x_1,...,x_\lambda )$ is an ideal of height $\lambda $ and $df(\lambda )=\lambda  $. Now let $p\geq 2$, since   $$ \ \ \Ical_{\Phi}=(x_{1}^{(p)},...,x_{m}^{(p)})\cap (x_{1}^{(p)},...,x_{m-1}^{(p)},\Ical_{\lambda_m})...
\cap (x_{1}^{(p)},...,x_{i-1}^{(p)},\Ical_{\lambda_i})\cap... (\Ical_{\lambda_m}),$$ we have that 
$$\hut \Ical_\Phi= \min \{\hut \Ical_{\lambda_i} + i-1 \},$$
by induction hypothesis $\hut \Ical_{\lambda_i}= df( \lambda_i)$ so $$\hut \Ical_\Phi= \min \{df( \lambda_i) + i-1 \}=df(\Phi ).$$ 
\item The proof is by induction on the number of generators  $ \mu (\Ical_\Phi )$ of the ideal $\Ical_\Phi$. 
The statement is clear if
$ \mu (\Ical_\Phi)=1.$ 

\noindent Suppose that $ \mu (\Ical_\Phi)>1.$ Let $\pi$  be a generator of $\Ical_\Phi$ being in  the last diagonal of $\Phi $, so we can write   $\pi =x_{i}^{(p)}g$ for some $i$, 
where $g\in \Ical_{\lambda_i }$ is in the last diagonal of $\lambda_i$.  By definition of a $p-$Ferrer tableau, 
 the ideal  generated by all the generators of $\Ical_\Phi$ except $x_{i}^{(p)}g$ is  a $p$-Ferrer ideal and 
we denoted it by $\Ical_{\Phi'}$. 

In the  example \ref{exo4322} we can perform several steps :
 
 $$\begin{array}{llllcllllcllllcllll}
4 & 3& 2& 2&&4 & 3& 2& 1&&4 & 3& 2& 1&&4 & 3& 2& 1  \\
3 & 2& 1& 0&&3 & 2& 1& 0&&3 & 2& 1& 0&&3 & 2& 1& 0   \\
2 & 0& 0& 0&\longrightarrow &2 & 0& 0& 0&\longrightarrow &2 & 0& 0& 0&\longrightarrow &2 & 0& 0& 0 \\
2 & 0& 0& 0&&2 & 0& 0& 0&&1 & 0& 0& 0&&0 & 0& 0& 0  \\
\\
&&\Phi&&&&&\Phi'\hskip -0.2cm&&&&&\Phi''\hskip -0.4cm&&&&&\hskip -0.4cm\Phi'''\\
\end{array}$$

let denote  $\alpha _{p}:=i$, so that 
 $$x_{i}^{(p)}g=x_{\alpha _{p}}^{(p)}x_{\alpha _{p-1}}^{(p-1)}...x_{\alpha _{1}}^{(1)}$$
 For any $k$ and $1\leq \beta <\alpha_k$ we have that $x_{\alpha _{p}}^{(p)}x_{\alpha _{p-1}}^{(p-1)}...
x_{\beta }^{(k)}...x_{\alpha _{1}}^{(1)}\in \Ical_{\Phi'}$, so we have that  
 $$(\{x_{1}^{(p)},...,x_{\alpha _{p}-1}^{(p)}  \},...,
\{x_{1}^{(1)},...,x_{\alpha _{1}-1}^{(1)}  \})\subset \Ical_{\Phi'}: x_{\alpha _{p}}^{(p)}...x_{\alpha _{1}}^{(1)}.$$

On the other hand let  $\Pi\in \Ical_{\Phi'}: x_{\alpha _{p}}^{(p)}...x_{\alpha _{1}}^{(1)}$  a monomial, we can suppose that no variable in $(\{x_{1}^{(p)},...,x_{\alpha _{p}-1}^{(p)}  \},...,
\{x_{1}^{(1)},...,x_{\alpha _{1}-1}^{(1)}  \})$ appears in $\Pi$, so $ \Pi x_{\alpha _{p}}^{(p)}...x_{\alpha _{1}}^{(1)}\in \Ical_{\Phi'}$ implies that there is a generator of $ \Ical_{\Phi'}$ of the type 
$x_{\beta  _{p}}^{(p)}...x_{\beta  _{1}}^{(1)} $ such that $\beta  _{i}\geq \alpha  _{i}$ for all $i=1,...,p$, this is in contradiction with the fact that $x_{\alpha _{p}}^{(p)}x_{\alpha _{p-1}}^{(p-1)}...x_{\alpha _{1}}^{(1)}$ is in the last diagonal of $\Phi $ and doesn't belongs to $\Ical_{\Phi'}$.
 In conclusion we  have that 
 $$\Ical_{\Phi'}: x_{\alpha _{p}}^{(p)}...x_{\alpha _{1}}^{(1)}=
(\{x_{1}^{(p)},...,x_{\alpha _{p}-1}^{(p)}  \},...,
\{x_{1}^{(1)},...,x_{\alpha _{1}-1}^{(1)}  \})$$
 is a linear ideal generated by   $\alpha _{p}+...+\alpha _{1}- (p)$ variables. Let remark that since 
$x_{\alpha _{p}}^{(p)}x_{\alpha _{p-1}}^{(p-1)}...x_{\alpha _{1}}^{(1)} $ is in the last diagonal the
 number of diagonals $\delta(\Phi ) $ in $\Phi $ is
$\alpha _{p}+...+\alpha _{1}- p+1$.  
 
We have the following exact sequence :
$$0\rightarrow S/(\Ical_{\Phi'}: (x_{i}^{(p)}g))(-p)\buildrel{\times x_{i}^{(p)}g}\over\longrightarrow S/(\Ical_{\Phi'})\rightarrow S/(\Ical_{\Phi})\rightarrow 0,$$
by applying the mapping cone construction we have that 
$$\beta_j(S/\Ical_{\Phi})= \beta_j(S/\Ical_{\Phi'})+ {{\delta(\Phi ) -1}\choose{j-1}}, \ \ \forall j=1,..., \projdim(S/\Ical_{\Phi}).$$ 
 By induction hypothesis the number of diagonals in $\Phi'$ 
 coincides with $\projdim(S/\Ical_{\Phi'})$. The number of diagonals in $\Phi'$ is either equal to the number of diagonals in $\Phi $ minus one,
 or
equal to the number of diagonals in $\Phi$. In both cases we have that 
$s_i(\Phi )=s_i(\Phi' )$ for $i=d-1,...,n-(\delta(\Phi )-1),$ 
 $s_{n-(\delta(\Phi ))}(\Phi' )=s_{n-(\delta(\Phi ))}(\Phi )-1$, and  $s_i(\Phi )=s_i(\Phi' )=0$ for $i<n-(\delta(\Phi ))$.
 
Let $c'=\hut \Ical_{\Phi'}$, It then follows that $$\beta_j(S/\Ical_{\Phi'})= 
{{c'+p-1}\choose{j+p-1}} {{j+p-2}\choose{p-1}}+\sum_{i=0}^{d-1}s_i(\Phi') {{n-i-1}\choose{j-1}},\ \ \forall j=1,..., \projdim(S/\Ical_{\Phi'}).$$
By induction hypothesis $\projdim(S/\Ical_{\Phi'})=\delta(\Phi' )$. We have to consider two cases:
\begin{enumerate}
\item $\delta (\Phi )=c$, this case can arrive only if the $c$ diagonal of $\Phi $ is full, so
$$c'=c-1,\ \  s_{n-c}(\Phi' )={{c-1+p-1}\choose{p-1}}-1, \delta(\Phi )=\delta(\Phi' )=c$$
$$\forall 1\leq j\leq c, \ \ \beta_j(S/\Ical_{\Phi})= \beta_j(c-1,p)+({{c-1+p-1}\choose{p-1}}-1){{c -1}\choose{j-1}}  + {{c -1}\choose{j-1}}.$$ 
$$\forall 1\leq j\leq c, \ \ \beta_j(S/\Ical_{\Phi})= \beta_j(c-1,p)+{{c-1+p-1}\choose{p-1}}{{c -1}\choose{j-1}} =\beta_j(c,p).$$
Let remark that by induction hypothesis $\beta_j(S/\Ical_{\Phi'})=0 $ for $j>c$, this implies that 
$\projdim(S/\Ical_{\Phi})=c=\delta(\Phi )$.
\item $\delta (\Phi )>c$, in this case $c'=c$
$$\beta_j(S/\Ical_{\Phi})= 
{{c+p-1}\choose{j+p-1}} {{j+p-2}\choose{p-1}}+\sum_{i=0}^{d-1}s_i(\Phi) {{n-i-1}\choose{j-1}}$$
and $\projdim(S/\Ical_{\Phi})= \projdim(S/\Ical_{\Phi'})$ equals the number of diagonals in $\Phi $.
\end{enumerate}
In particular it follows that if  the number of diagonals in $\Phi'$ is  equal to the number of diagonals in $\Phi $
 minus one, $\projdim(S/\Ical_{\Phi})= \projdim(S/\Ical_{\Phi'})+1$ is the number of diagonals in $\Phi $. If 
 the number of diagonals in $\Phi'$ is  equal to the number of diagonals in $\Phi $, then 
 $\projdim(S/\Ical_{\Phi})= \projdim(S/\Ical_{\Phi'})$ equals the number of diagonals in $\Phi $.
\begin{proposition} $\ara(\Ical_\Phi )=\cd(\Ical_\Phi )=\pd (S/\Ical_\Phi )$. 

\end{proposition}
\demo Recall that a monomial in the $p-$Ferrer ideal (or tableau) 
$x_{\alpha _{p}}^{(p)}x_{\alpha _{p-1}}^{(p-1)}...x_{\alpha _{1}}^{(1)} $ is in the 
$\alpha _{p}+\alpha _{p-1}+...+\alpha _{1}-p+1$ diagonal.
Let $\Kcal_j $ the set of all monomials in the Ferrer tableau lying in the $j$ diagonal
and let $F_j=\sum_{M\in \Kcal_j }M$, we will prove that for any $M\in \Kcal_j$ , we have$M^2\in (F_1,...,F_j)$.
If $j= 2$ , let $
M= x_{\alpha _{p}}^{(p)}x_{\alpha _{p-1}}^{(p-1)}...x_{\alpha _{1}}^{(1)} $, with 
$\alpha _{p}+\alpha _{p-1}+...+\alpha _{1}-p+1=2$, then 
$$MF_2= M^2+\sum (x_{\alpha _{p}}^{(p)}x_{\alpha _{p-1}}^{(p-1)}...x_{\alpha _{1}}^{(1)})M'$$
One monomial  $M' \in \Kcal_2$, $M'\not= M$ can be written 
$$x_{\beta  _{p}}^{(p)}x_{\beta  _{p-1}}^{(p-1)}...x_{\beta_1 ^{(1)}}$$
with $\beta  _{p}+\beta  _{p-1}+...+\beta _{1}-p+1=2$, this implies that
$\beta  _{i}=1$ for all $i$ except one value $i_0$, for which $\beta  _{i_0}=2$  and also 
$\alpha  _{j}=1$ for all $j$ except one value $j_0$, for which $\alpha  _{j_0}=2$. Since 
$M'\not= M$ we must have $x_{1}^{(P)}x_{1}^{(p-1)}...x_{1}^{(1)} $ divides $MM'$.

Now let $j\geq 3$, let 
$M= x_{\alpha _{p}}^{(p)}x_{\alpha _{p-1}}^{(p-1)}...x_{\alpha _{1}}^{(1)} $, with 
$\alpha _{p}+\alpha _{p-1}+...+\alpha _{1}-p+1=j$, then 
$$MF_2= M^2+\sum (x_{\alpha _{p}}^{(p)}x_{\alpha _{p-1}}^{(p-1)}...x_{\alpha _{1}}^{(1)})M'$$
One monomial $M' \in \Kcal_j$, $M'\not= M$  
can be written $$x_{\beta  _{p}}^{(p)}x_{\beta  _{p-1}}^{(p-1)}...x_{\beta_1^{(1)}}$$
with $\beta  _{p}+\beta  _{p-1}+...+\beta _{1}-p+1=j$, let $i_0$ such that $\beta _{i_0}\not=\alpha _{i_0}$
 if
$\beta  _{i_0}<\alpha _{i_0}$ then $\displaystyle\frac{M}{x_{\alpha _{i_0}}^{(i_0)}}{x_{\beta _{i_0}}^{(i_0)}}\in K_i$
 for some $i<j$, and if 
 $\beta  _{i_0}>\alpha _{i_0}$ then $\displaystyle\frac{M'}{x_{\beta _{i_0}}^{(i_0)}}{x_{\alpha _{i_0}}^{(i_0)}}
\in K_i$
 for some $i<j$,in both cases $MM'\in (K_i)$ for some $i<j$.
 As a consequence $\ara (\Ical_{\Phi})\leq \pd(S/\Ical_{\Phi})$, but $\Ical_{\Phi}$ is a monomial ideal, so by a Theorem
 of Lyubeznik $\cd (\Ical_{\Phi})= \pd(S/\Ical_{\Phi})$, and $\cd (\Ical_{\Phi})\leq \ara (\Ical_{\Phi})$,
 so we have the equality $\ara (\Ical_{\Phi})= \pd(S/\Ical_{\Phi})$. Let remark that the  equality 
$\cd (\Ical_{\Phi})= \pd(S/\Ical_{\Phi})$ can be recovered by direct computations in the case of $p-$Ferrer ideals.
 \end{enumerate}
 The reader should consider the relation between our theorem and the  following result from \cite{eg}:
 \begin{proposition} If $R:=S/\Ical$ is a homogeneous ring with $p-$linear resolution over an infinite field, and $x_i\in R_1$ are elements such that $x_{i+1}$ is a non zero divisor on $(R/(x_1,...,x_i))/H^0_{\bf m}(R/(x_1,...,x_i),$ where ${\bf m}$ is the unique homogeneous maximal ideal of $S$, then 
\begin{enumerate}
\item $s_i(R)=length(H^0_{\bf m}(R/(x_1,...,x_i)_{p-1}),$ for $i=0,...,\dim R-1$.
\item If $R$ is of codimension $c$, and $n:=\dim S$, the Betti numbers of $R$ are given by:
 $${\rm for\ \ }j =1,..., n-\depth(R)\ \ \beta_j(R)= \beta_j(c,p)
+\sum_{i=0}^{d-1}s_i {{n-i-1}\choose{j-1}},$$
\end{enumerate}
where $\beta_j(c,p)={{c+p-1}\choose{j+p-1}} {{j+p-2}\choose{p-1}}$ are the betti numbers of a Cohen-Macaulay ring having $p-$linear resolution, of codimension $c.$ 
\end{proposition}
 We have the following corollary:
\begin{corollary} If $R:=S/\Ical$ is a homogeneous ring with $p-$linear resolution over an infinite field,of codimension $c,$ and $n:=\dim S$, then 
$$\beta_j(c,p)\leq \beta_j(R)\leq \beta_j(n-\depth(R),p).$$
\end{corollary}
\demo As a consequence of the above proposition we have that
 $s_i\leq \displaystyle{{n-(i+1) +p-1}\choose{p-1}}$ so that
 $$\beta_j(c,p)\leq \beta_j(R)\leq  \beta_j(c,p)
+\sum_{i=0}^{d-1} \displaystyle{{n-(i+1) +p-1}\choose{p-1}}{{n-i-1}\choose{j-1}}$$
 By direct computations we have that $\beta_j(c,p)+\displaystyle{{n-d +p-1}\choose{p-1}}\displaystyle{{n-d}\choose{j\ell -1}} =\beta_j(c+1,p)$, which implies
 $$\beta_j(c,p)\leq \beta_j(R)\leq  \beta_j(c+1,p)
+\sum_{i=0}^{d-2} \displaystyle{{n-(i+1) +p-1}\choose{p-1}}{{n-i-1}\choose{j-1}},$$
by repeating the above computations we got the corollary.
 

 \section{Hilbert series of ideals with $p-$linear resolution.}
 Let $\Ical\subset S$ be an ideal with $p-$linear resolution, it follows from \cite{eg}, that the Hilbert series of $S/\Ical$ is given by
 $$H_{S/\Ical} (t)= \displaystyle\frac{\displaystyle\sum_{i=0}^{p-1}{{c+i-1}\choose{i}}t^i-t^p\Big(\displaystyle\sum_{i=1}^{ d} s_{d-i} (1-t)^{i-1}\Big)}{(1-t)^{d} }$$
where $d=n-c$
 In the case where the ring $S/\Ical$ is Cohen-Macaulay, we have :
 $$H_{S/\Ical} (t)= \displaystyle\frac{\displaystyle\sum_{i=0}^{p-1}{{c+i-1}\choose{i}}t^i}{(1-t)^{d} }$$
 
\begin{definition}For any non zero natural numbers $c,p$, we  set $$h(c,p)(t):=\displaystyle\sum_{i=0}^{p-1}{{c+i-1}\choose{i}}t^i.$$
\end{definition}
Remark that  the $h-$vector of the polynomial $h(c,p)(t)$ is log concave, since for $i=0,...,p-3$, we have that
 $${{c+i-1}\choose{i}} {{c+i+1}\choose{i+2}}\leq ({{c+i}\choose{i+1}})^2.$$

\begin{lemma} For any non zero natural numbers $c,p$, we have the relation 
$$ 1-h(c,p)(1-t)t^c=h(p,c)(t)(1-t)^p , $$
in particular $h(c,p)(t)\ (1-t)^c=1-h(p,c)(1-t)\ t^p$, $h(c,p)(t)\ (1-t)^c\equiv 1 \mod t^p$.
 
\end{lemma}
\demo Let $\Ical $ be a square free monomial ideal  having a $p-$linear resolution, such that $S/\Ical$ is a Cohen-Macaulay ring of codimension $c$, let   $\Jcal:=\Ical^*$ be the Alexander dual of $\Ical $, it then follows that $S/\Jcal$ is a Cohen-Macaulay ring of codimension  $p$ which has    a $c-$linear resolution.
$$\displaystyle H_{S/\Ical} (t)=\frac{h(c,p)(t)}{(1-t)^{n-c}}=\displaystyle\frac{1-B_{S/\Ical}(t)}{(1-t)^{n}}$$
$$\displaystyle H_{S/\Jcal} (t)=\frac{h(p,c)(t)}{(1-t)^{n-p}}=\displaystyle\frac{1-B_{S/\Jcal}(t)}{(1-t)^{n}}$$
and by Alexander duality on the Hilbert series we have that :
$1-B_{S/\Ical}(t)=B_{S/\Jcal}(1-t)$ but $h(c,p)(t)(1-t)^{c}=1-B_{S/\Ical}(t)$ and
 $h(p,c)(t)(1-t)^{p}=1-B_{S/\Jcal}(t)$, so $B_{S/\Jcal}(1-t)=1-h(p,c)(1-t)(t)^{p}$, so our claim follows from these identities.

\begin{corollary}Let $\Ical\subset S$ be any homogeneous ideal, $c=\hut (\Ical) , d=n-c$ and $p$ the smallest degree of a set of generators. Then we can write  $H_{S/\Ical} (t)$ as follows 
$$H_{S/\Ical} (t)= \displaystyle\frac{h(c,p)(t)-t^p\Big(\displaystyle\sum_{i=1}^{ \delta(\Ical) } s_{\delta(\Ical)-i} (1-t)^{i-1}\Big)}{(1-t)^{d} },$$
where the numbers $s_0,...,s_{\delta(\Ical)-1}$ are uniquely determined.
\begin{enumerate}
\item Let $\Jcal $ be a square free monomial ideal such that $S/\Jcal$ is a Cohen-Macaulay ring of codimension $p$, let   $\Ical:=\Jcal^*$ be the Alexander dual of $\Jcal $, it then follows that $S/\Ical$ has  a $p-$linear resolution. Let $c=\codim (S/\Ical)$.
Then 
$$H_{S/\Ical} (t)=\displaystyle\frac{ h(c,p)(t)-t^p\Big(\displaystyle\sum_{i=1}^{ d} s_{d-i} (1-t)^{i-1}\Big)}{(1-t)^{n-c}}, $$ $$ H_{S/\Jcal} (t)=\displaystyle\frac{h(p,c)(t)+t^c\Big(\displaystyle\sum_{i=1}^{ d} s_{d-i} t^{i-1}\Big)}{(1-t)^{n-p}}$$

\item Let $\Ical $ be any square free monomial ideal $c=\codim (S/\Ical)$, $p$ the smallest degree of a set of generators of $\Ical $.  Let   $\Jcal:=\Ical^*$ be the Alexander dual of $\Ical $, then $p=\codim (S/\Jcal)$,  $c$ is the smallest degree of a set of generators of $\Jcal $ and 
$$H_{S/\Ical} (t)=\displaystyle\frac{ h(c,p)(t)-t^p\Big(\displaystyle\sum_{i=1}^{ \delta(\Ical) } s_{\delta(\Ical)-i} (1-t)^{i-1}\Big)}{(1-t)^{n-c}}, $$ $$ H_{S/\Jcal} (t)=\displaystyle\frac{h(p,c)(t)+t^c\Big(\displaystyle\sum_{i=1}^{ \delta(\Ical)} s_{\delta(\Ical)-i} t^{i-1}\Big)}{(1-t)^{n-p}}.$$
\end{enumerate}

\end{corollary}
\demo Since $1,t,..., t^{p},t^{p}(1-t),...,t^{p}(1-t)^{k},...,$
 are linearly independent the numbers $s_i$ are uniquely defined.

$$\displaystyle H_{S/\Ical} (t)=\frac{h_{S/\Ical} (t)}{(1-t)^{n-c}}=\displaystyle\frac{1-B_{S/\Ical}(t)}{(1-t)^{n}}$$
$$\displaystyle H_{S/\Jcal} (t)=\frac{h_{S/\Jcal}(t)}{(1-t)^{n-p}}=\displaystyle\frac{1-B_{S/\Jcal}(t)}{(1-t)^{n}}$$
 by Alexander duality on the Hilbert series we have that :
$$B_{S/\Jcal}(t)=1-B_{S/\Ical}(1-t)=(h(c,p)(1-t)-(1-t)^p\Big(\displaystyle\sum_{i=1}^{ \delta(\Ical) } s_{\delta(\Ical)-i} t^{i-1}\Big))t^c,$$	  but $ h(c,p)(1-t)t^c=1-h(p,c)(t)\ (1-t)^p,$ so 
$$1-B_{S/\Jcal}(t)= (h(p,c)(t)+ t^c\Big(\displaystyle\sum_{i=1}^{ \delta(\Ical) } s_{\delta(\Ical)-i} t^{i-1}\Big) )(1-t)^p. $$ This proves the claim.

\begin{theorem}
\begin{enumerate}
\item For any $M-$vector ${\bf h}=(1,h_1,...)$  there exists  $\Phi $   a $p-$Ferrer tableau such that $h_i$ counts
the number of elements in the $i-$diagonal of $\Phi $.

\item the $h-$vector of any $p-$regular ideal is the $h-$vector of a $p-$Ferrer ideal.
\item For any $M-$vector ${\bf h}=(1,h_1,...)$  we can explicitely construct a $p-$Ferrer tableau  $\Phi $ such that
${\bf h}=(1,h_1,...)$ is the $h-$vector of $\Ical_\Phi^*$. 
\end{enumerate}

\end{theorem}
 \demo \begin{enumerate}
 \item Let ${\bf h}=(1,h_1,...)$ be the $h-$vector of $S/\Jcal$, by Macaulay, \cite{st} 2.2 theorem $h$ is obtained as the $M-$vector of a multicomplex $\Gamma $, where $h_i$ counts the monomials of degree $i$ in $\Gamma $.
We establish  a correspondence between  multicomplex $\Gamma $ and $p-$ Ferrer ideals:

Suppose that $\Gamma $ is a multicomplex in the variables $x_1,...,x_n$, to any monomial $x_1^{\alpha _1}...x_n^{\alpha _n}\in \Gamma $ we associated the vector $(\alpha _1+1,...,\alpha _n+1)\in (\bn^*)^n$, 
let $\Phi $ be the image of $\Gamma $. By definition $\Gamma $ is a multicomplex if and only if for any $u\in \Gamma$, and if  $v$ divides $u$ then $v\in \Gamma$, this property is equivalent to the property: 

For any $(\alpha _1+1,...,\alpha _n+1)\in \Phi $ and $(\beta _1+1,...,\beta _n+1) \in (\bn^*)^n$ such that $\beta _i\leq \alpha_i$ for all $i$ then  $(\beta _1+1,...,\beta _n+1)\in \Phi . $ That is  $\Phi $ is a $p-$Ferrer tableau,   such that $h_i$ counts the number of elements in the $i-$diagonal of $\Phi $.
\item 
 Let $\Ical\subset S$ be any graded ideal with
 $p-$linear resolution, 
let $Gin(\Ical)$ be the generic initial, by a theorem of Bayer and Stillman, $Gin(\Ical)$ has a $p-$linear resolution, on the other hand they have the same Hilbert series, and from the remark in the introduction they have the same betti numbers. $Gin(\Ical)$  is a monomial ideal, we can take the polarisation $P(Gin(\Ical))$, this is a square free monomial having $p-$linear resolution and the same betti numbers as $Gin(\Ical)$,
 the Alexander dual $P(Gin(\Ical))^*$ is Cohen-Macaulay of codimension $p$
, so there exists a Ferrer tableau  $\Phi$ such that the $h-$vector of $S/P(Gin(\Ical))^*$ 
is the generating function of the diagonals of $\Phi $, moreover  the $h-$vector of $S/P(Gin(\Ical))^*$coincides with the $h-$vector of $S/(\Ical_\Phi )^*$. By the above proposition the $h-$vector of $S/P(Gin(\Ical))^*$ determines uniquely the  $h-$vector of $S/P(Gin(\Ical))$, and the last one coincides with the $h-$vector of $S/\Ical_\Phi $.
\item Let recall from \cite{st} how to associate to a $M-$vector ${\bf h}=(1,h_1,...,h_l)$ a multicomplex 
$\Gamma_{\bf h}$. 
For all $i\geq 0$ list all monomials in $h_1$ variables in reverse lexicographic order, let 
$\Gamma_{{\bf h},i}$ be set of first $h_i$ monomials in this order, and $\Gamma_ {\bf h}=\bigcup_{i=0}^{i=l} \Gamma_ {{\bf h},i}, $ in the first item we have associated to a multicomplex a $p-$Ferrer tableau $\Phi$ such that $h_i$ is the number of elements in the $i-$diagonal of $\Phi$. By the second item the $h-$vector of $S/(\Ical_\Phi )^*$  is exactly ${\bf h}$.  
\end{enumerate}  
\begin{example} We consider the $h-$vector, $(1,4,3,4,1)$, following \cite{st}, this $h-$vector corresponds to the
multicomplex 
$$1;x_1,...,x_4;x_1^2,x_1x_2,x_2^2;
x_1^3,x_1^2x_2,x_1x_2^2,x_3^2;x_1^4,$$ 
and to the following $p-$Ferrer ideal $\Ical_\Phi$ generated by:
$$ s_1t_1u_1v_1,$$ 
$$s_2t_1u_1v_1,s_1t_2u_1v_1,s_1t_1u_2v_1,s_1t_1u_1v_2,$$
$$s_3t_1u_1v_1,s_2t_2u_1v_1,s_1t_3u_1v_1,$$
$$s_4t_1u_1v_1,s_3t_2u_1v_1,s_2t_3u_1v_1,s_1t_4u_1v_1,$$
$$s_5t_1u_1v_1,$$
$\Ical_\Phi$ has the following prime decomposition:
$$(v_1,v_2)\cap (u_1,v_1)\cap (s_1,v_1)\cap (t_1,v_1)\cap (u_1,u_2)\cap (s_1,u_1)\cap (t_1,u_1)\cap $$
$$\cap (t_1,t_2,t_3,t_4)\cap (t_1,t_2t_3,s_1)\cap (t_1,t_2,s_1,s_2)\cap (s_1,s_2,s_3,s_4,s_5)$$
and $\Ical_\Phi^*$ is generated by
$$ v_1 v_2 ,  u_1 v_1 ,  s_1 v_1 ,  t_1 v_1 ,  u_1 u_2 ,  s_1 u_1 ,  t_1 u_1 , $$
$$  t_1 t_2t_3 t_4 ,  t_1 t_2t_3 s_1 ,  t_1 t_2 s_1 s_2,  s_1 s_2 s_3 s_4 s_5 $$
and the $h-$vector of $S/(\Ical_\Phi )^*$ is $(1,4,3,4,1)$.
\end{example}
\section{Examples }
Let
$S=K[x_1,...,x_n]$ be a polynomial ring. Let $\alpha \in \bn^*$, for any element $P\in S$ we set $\widetilde P(x)=P(x_1^\alpha,...,x_n^\alpha )$, more generally for any matrix with entries in $S$ we set $\widetilde M$ be matrix obtained by changing the entry $P_{i,j}$ of $M$ to $ \widetilde P_{i,j} $.
\begin{lemma}Suppose that $$F^{\bullet  }\  :\ 0\rightarrow F_s\buildrel{M_s}\over\rightarrow F_{s-1}\rightarrow ....\rightarrow F_1\buildrel{M_1}\over\rightarrow  F_0\rightarrow 0$$ 
is a minimal free resolution of a graded $S-$module $M$, then 
$$\widetilde  F^\bullet  :0\rightarrow \widetilde  F_s\buildrel{\widetilde M_s}\over\rightarrow \widetilde  F_{s-1}\rightarrow ....\rightarrow \widetilde  F_1\buildrel{\widetilde M_1}\over\rightarrow  \widetilde  F_0\rightarrow 0$$ is a minimal free resolution of a graded $S-$module $\widetilde  M$. 
If $F^{\bullet  }$ is a pure free resolution, that is  $F_{i}=S^{\beta _i}(-a_i)$ for all $i=0,...,s$, then $\widetilde  F^\bullet $ is also pure and  $\widetilde  F_{i}=S^{\beta _i}(-a_i\alpha )$ for all $i=0,...,s$. 
\end{lemma}
\begin{corollary}\begin{enumerate}
\item Let $\Phi(p,c) $ be the  $p-$ Ferrer diagram Cohen-Macaulay of codimension $c$. Let $\widetilde \Phi(p,c) $ be the Ferrer $p-$ Ferrer diagram obtained from $\Phi(p,c) $ by dividing any length unit into $\alpha $ parts, then  the Alexander dual $\Ical_{\widetilde\Phi(p,c)}^*$ has a pure resolution of type $(0,c\alpha,...,(c+p-1)\alpha )$.  
\item Let consider any sequence $0<a_1<a_2$, suppose that $\beta_0-\beta _1t^{a_1} +\beta _2t^{a_2}$ is the Betti polynomial of a Cohen-Macaulay module of codimension 2, then we must have:
 $$ a_1= c\alpha   , a_2=(c+1)\alpha  , \beta _1=\displaystyle\frac{a_2}{a_2-a_1}\beta_0, \beta _2=\displaystyle\frac{a_1}{a_2-a_1}\beta_0.$$
if $a_2-a_1$ is a factor of $ a_2,a_1$ then we can write $$ a_1= c\alpha   , a_2=(c+1)\alpha  , \beta _1=(c+1)\beta_0, \beta _2=c\beta_0,$$ with $c$ a natural number.
In particular a module obtaining by taking $\beta _0$ copies of $\Ical_{\widetilde\Phi(2,c)}^*$, has a pure resolution of type $(0,c\alpha,(c+1)\alpha )$. 
\end{enumerate}

\end{corollary}
\begin{example} Let $S=K[a,b,c]$, consider the free resolution of the algebra $S/(ab,ac,cd)$:  $$0\longrightarrow S^{2}\buildrel{\pmatrix{a&0\cr -d&b\cr 0&-c\cr } }\over\longrightarrow S^3\buildrel{\pmatrix{cd&ac&ab\cr}}\over\longrightarrow S\longrightarrow 0$$ then we have a pure free resolution
$$0\longrightarrow S^{2\beta _0}\buildrel{M_1 }\over\longrightarrow S^{3\beta _0}\buildrel{M_0}\over\longrightarrow S^{\beta _0}\longrightarrow 0,$$
where $$M_1=\pmatrix{a^\alpha &0&&&&&\cr -d^\alpha &b^\alpha &&&&&\cr 0&-c^\alpha &&&&&\cr&&.&&&&\cr&&&.&&&\cr&&&&.&&\cr&&&&&a&0\cr &&&&&-d^\alpha &b^\alpha \cr&&&&&0&-c^\alpha \cr }, M_0=\pmatrix{c^\alpha d^\alpha &a^\alpha c^\alpha &a^\alpha b^\alpha &&&&&&\cr&&&.&&&&\cr&&&&.&&&\cr&&&&&.&&\cr&&&&&&cd&a^\alpha c^\alpha &a^\alpha b^\alpha \cr },$$ with the obvious notation.
\end{example}
\begin{example}
 The  algebra $S/(ab,ae,cd,ce,ef)$ has Betti-polynomial $1-5t^2+5t^3-t^5$ but has not pure resolution.

\end{example}

\begin{example}{\bf Magic squares} Let $S$ be a polynomial ring of dimension  $n!$, It follows from \cite{st} that the toric ring of $n\times n$ magic squares is a quotient $R_{\Phi_n }=S/\Ical_{\Phi_n }$,  its  $h-$polynomial is as follows:
$$h_{R_{\Phi_n }}(t)= 1+h_1t+...+h_lt^d,$$
where $h_1=n!-(n-1)^2-1, d=(n-1)(n-2) $,
If $n=3$ we have $h_1=1, d=2$, and there is no relation of degree two between two permutation matrices, but we have a degree three relation. Set
$$M_1=\pmatrix{1&0&0\cr 0&1&0\cr0&0&1\cr}, \ \ M_2=\pmatrix{0&0&1\cr 1&0&0\cr0&1&0\cr},\ \ M_3=\pmatrix{0&1&0\cr 0&0&1\cr 1&0&0\cr},$$
$$M_4=\pmatrix{1&0&0\cr 0&0&1\cr0&1&0\cr}, \ \ M_5=\pmatrix{0&1&0\cr 1&0&0\cr0&0&1\cr},\ \ M_6=\pmatrix{0&0&1\cr 0&1&0\cr 1&0&0\cr},$$
We can see that $M_1+M_2+M_3=M_4+M_5+M_6$, so this relation gives a degree three generator in $\Ical_{\Phi_3 }$, and in fact $\Ical_{\Phi_3 }$ is generated by this relation. 
By using the cubic generators of $\Ical_{\Phi_3 }$ we get cubic generators of $\Ical_{\Phi_n }$ for  $n\geq 4$, but we have also quadratic generators, for example :
$$\pmatrix{1&0&0&0\cr 0&1&0&0\cr 0&0&1&0\cr 0&0&0&1\cr}+ \pmatrix{0&0&1&0\cr 0&0&0&1\cr 1&0&0&0\cr 0&1&0&0\cr}=
\pmatrix{1&0&0&0\cr 0&0&0&1\cr 0&0&1&0\cr 0&1&0&0\cr}+ \pmatrix{0&0&1&0\cr 0&1&0&0\cr 1&0&0&0\cr 0&0&0&1\cr}
$$
So  for  $n\geq 4$,  the smallest degree of a generator of the toric ideal  $\Ical_{\Phi_n }$ is of degree 2.
and unfortunately our proposition can  give   only information about $h_1$.
\end{example}

\begin{example} 
The Hilbert series of the following $p-$Ferrer tableaux are respectively: $$\displaystyle\frac{1+3t+6t^2}{(1-t)^6},\ \  \displaystyle\frac{1+2t+3t^2-5t^3}{(1-t)^7},\ \ \displaystyle\frac{1+3t+6t^2-t^3}{(1-t)^6}.$$
Let remark that $\displaystyle\frac{1+2t+3t^2-6t^3}{(1-t)^7}= \frac{1+3t+6t^2}{(1-t)^6} .$
  \begin{center}
\includegraphics[height=1.6 in,width=1.6in]{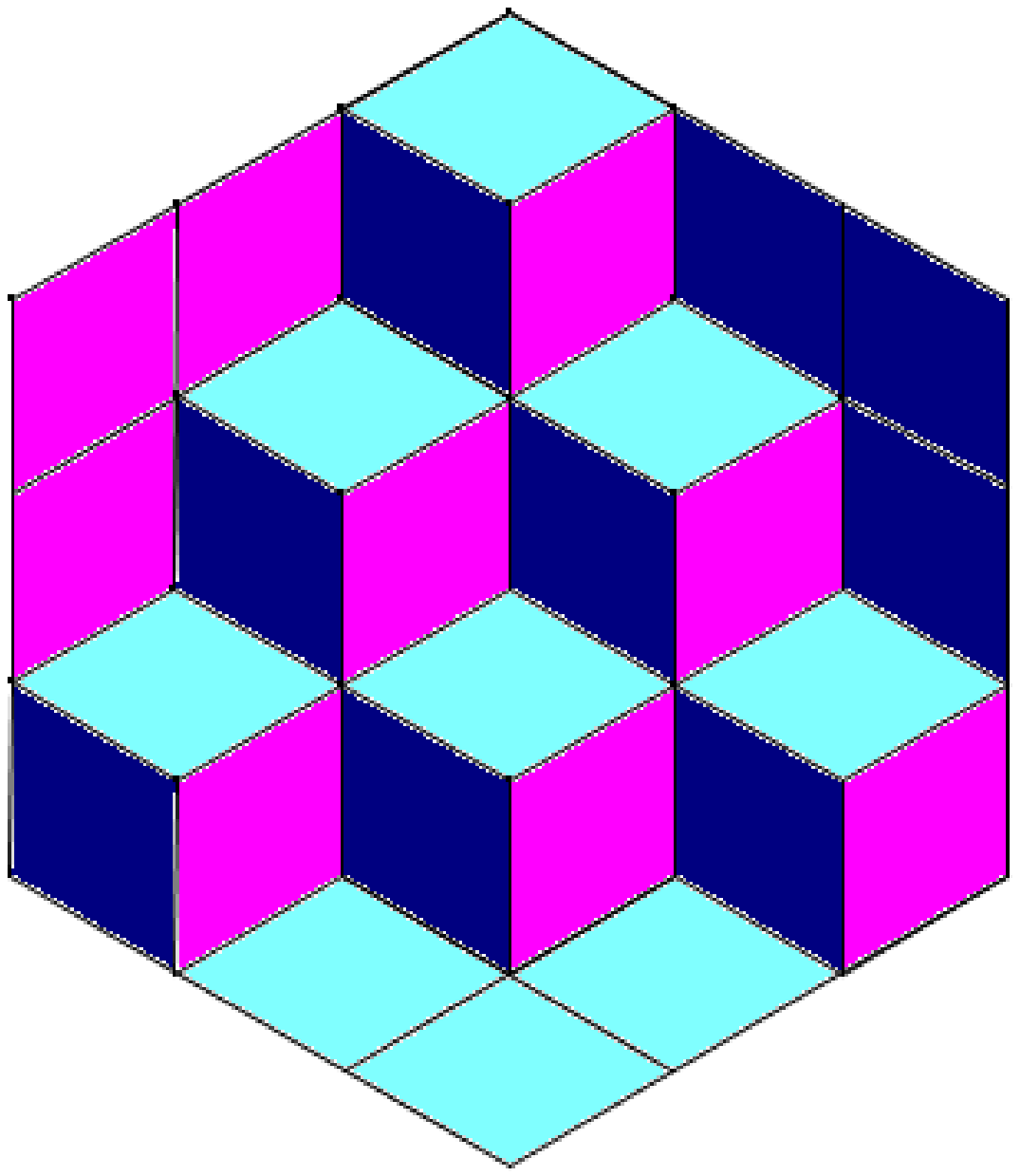}
\includegraphics[height=1.6 in,width=1.6in]{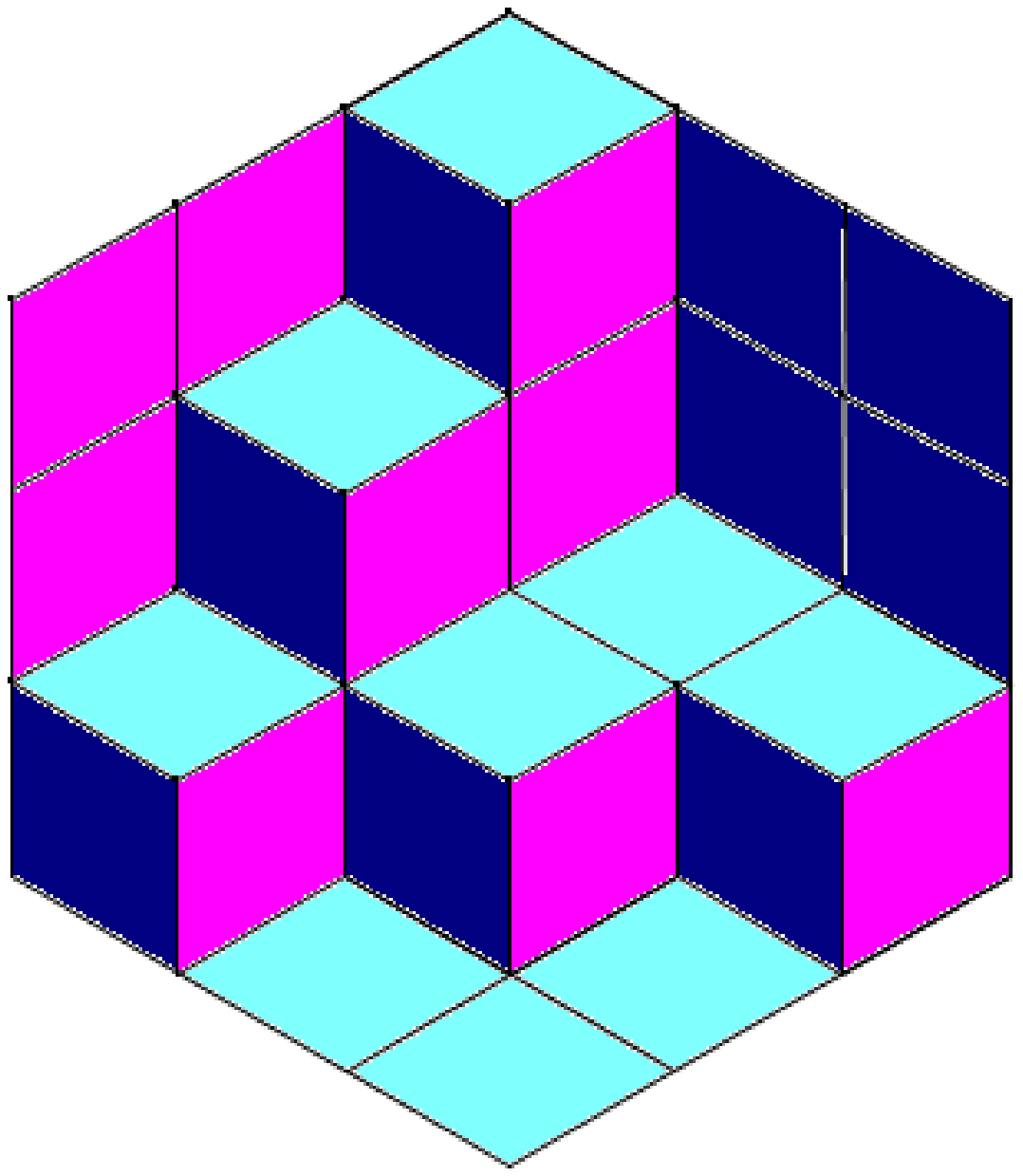}
\includegraphics[height=1.5 in,width=1.5in]{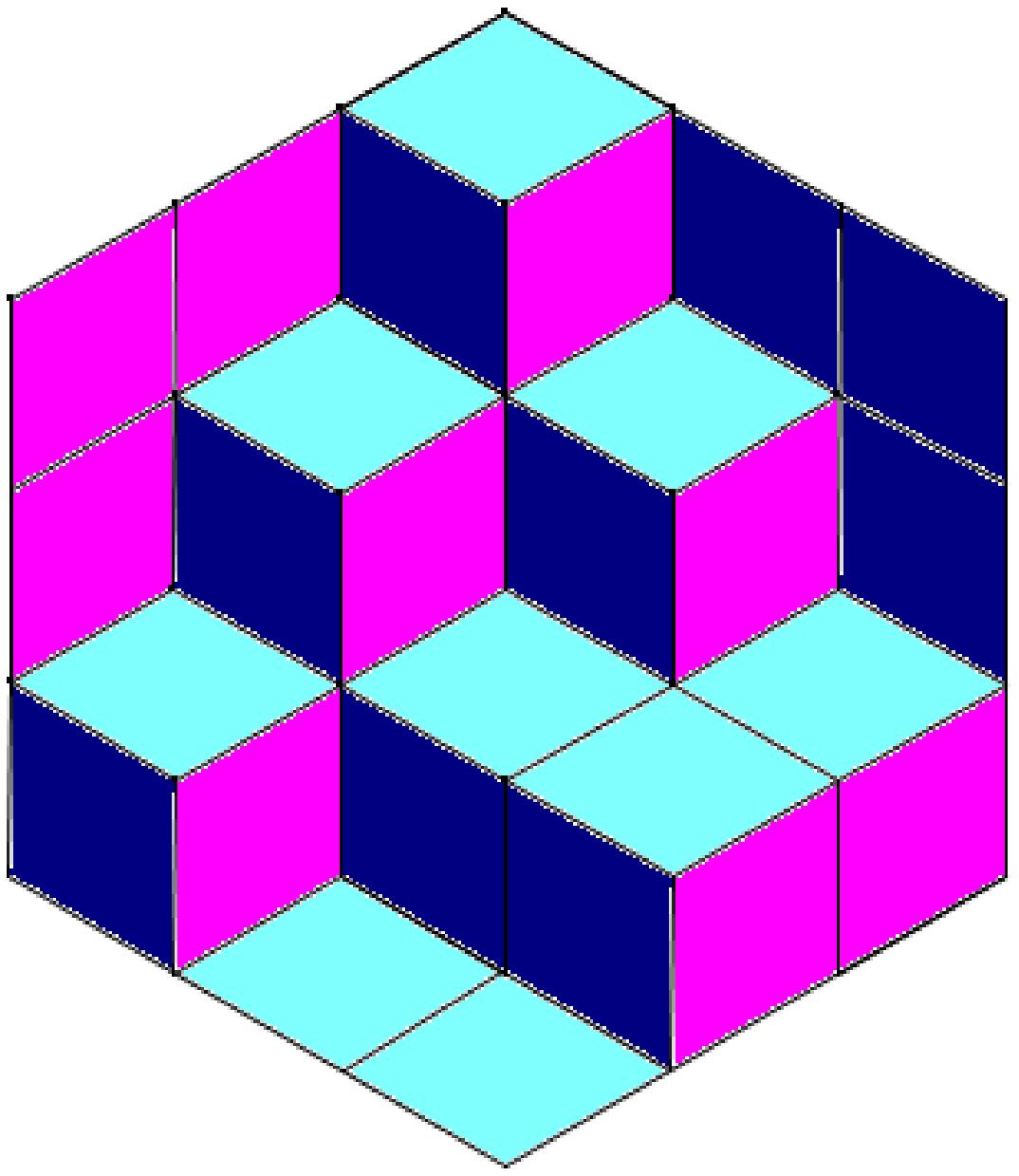}
\end{center}

\end{example} 
\begin{example} The Hilbert series of the following $p-$Ferrer tableaux are respectively: $$  H(t)=\displaystyle\frac{1+t+t^2-t^3(2+ 2(1-t))}{(1-t)^6},
 H(t)=\displaystyle\frac{1+t+t^2-t^3(2+ 3(1-t)+(1-t)^2)}{(1-t)^6}$$
\begin{center}
\includegraphics[height=1.6 in,width=1.6in]{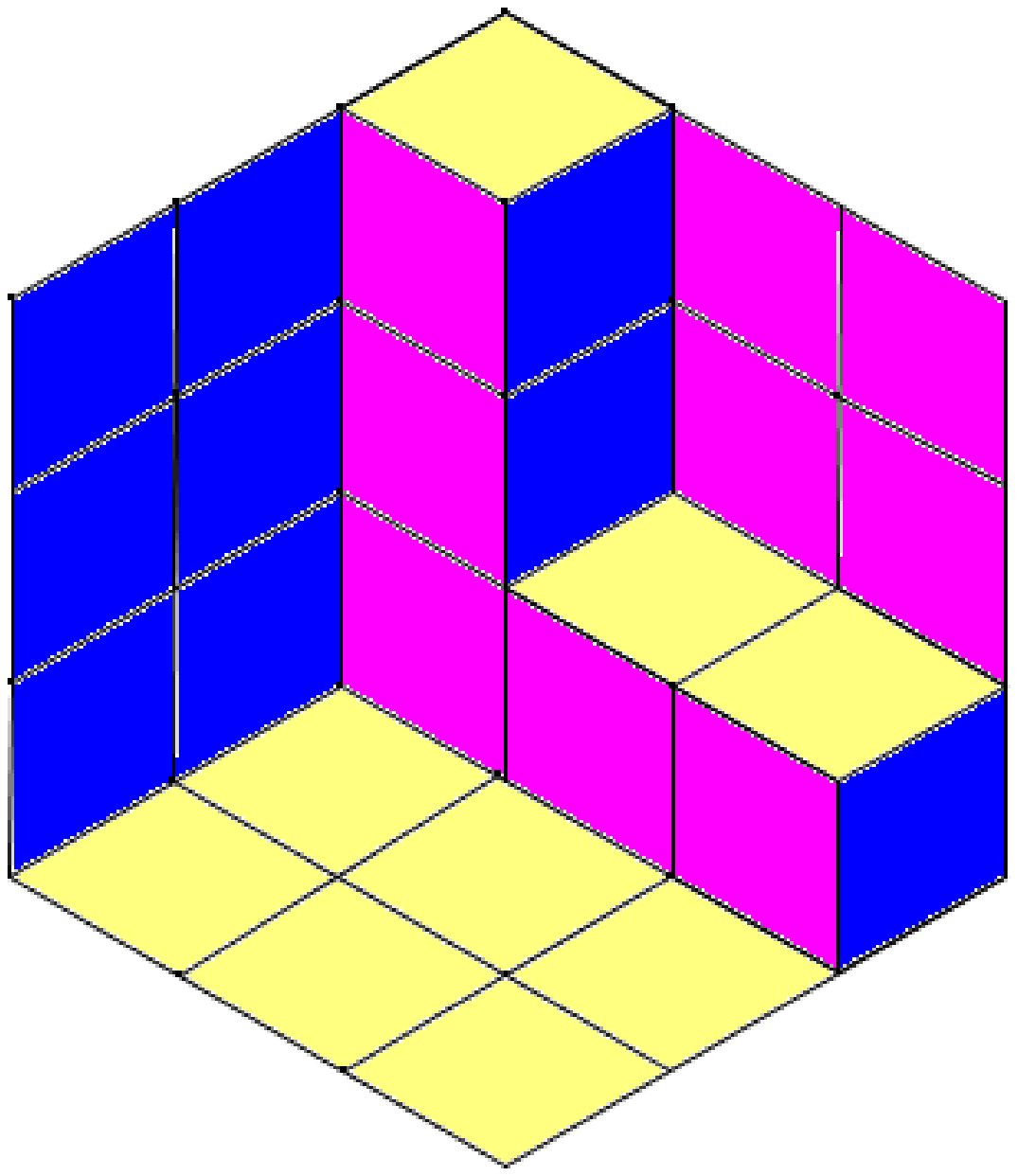}
\includegraphics[height=1.6 in,width=1.6in]{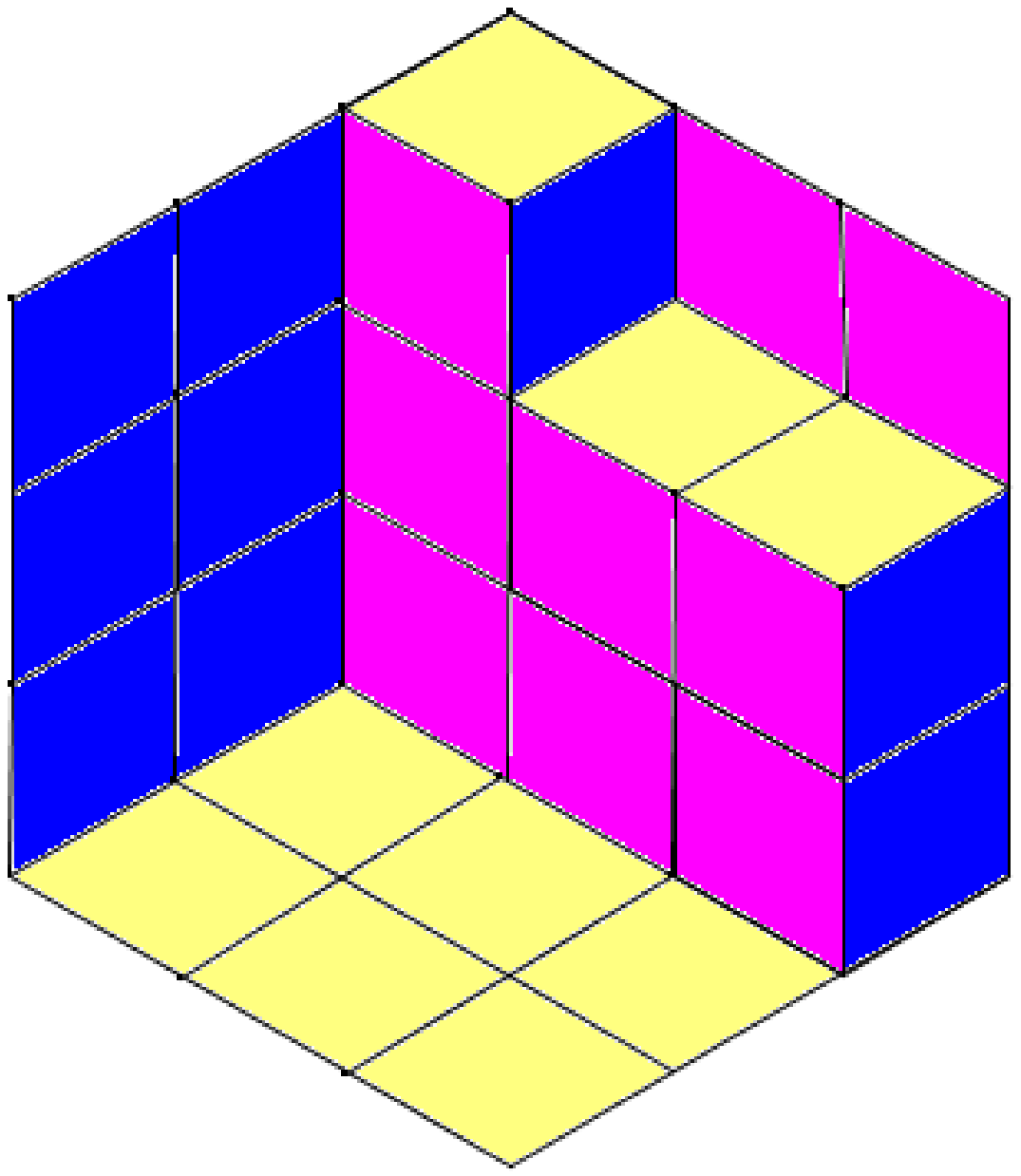}
\end{center}

\end{example}

 \end{document}